\documentclass{article}
\usepackage{amsmath}
\usepackage{amsfonts}
\usepackage{amssymb}
\usepackage{pstricks,pst-node,amssymb,amsmath,graphics,latexsym,tabularx,shapepar}
\usepackage{graphicx}
\newcommand{\bpartial}{\mathop{\partial\kern -4pt\raisebox{.8pt}{$|$}}}
\newcommand{\bra}{\mathopen{[\kern-1.6pt[}}
\newcommand{\ket}{\mathclose{]\kern-1.5pt]}}
\newcommand{\bbra}{\mathopen{[\kern-2.2pt[\kern-2.3pt[}}
\newcommand{\bket}{\mathclose{]\kern-2.1pt]\kern-2.3pt]}}

\makeindex

\title{Geometric Dynamical  Systems}
\author{Ghorbanali Haghighatdoost\\
\small Azarbaijan Shahid Madani University, Tabriz, Iran\\
\small Email: gorbanali@azaruniv.ac.ir}
\date{}

\begin{document}
\maketitle

\begin{abstract}
This article provides a conceptual and historical review of the evolution of integrable Hamiltonian systems from the Moscow School of A.~T.~Fomenko to the emerging Azarbaijan School of Geometric Dynamical Systems founded by the author. Beginning with the topological classification of integrable systems through Liouville foliations, atoms, and molecular invariants, the paper traces how these geometric ideas evolved into modern frameworks based on Lie groupoids, Lie algebroids, and fractional calculus. The author’s doctoral dissertation at Moscow State University (2004) extended the Fomenko theory to new integrable systems on $\mathfrak{so}(4)$, constructing novel topological molecules and describing the hierarchy of singularities and bifurcations. Upon his return to Iran, he established a comprehensive research program at Azarbaijan Shahid Madani University, integrating topology, geometry, and analysis to form a coherent Iranian branch of the global theory of integrable systems.
This program unifies the classical and the modern: from the Euler–Poincaré and Hamilton–Jacobi formalisms on Lie groupoids and algebroids, to fractional and quantum mechanical models involving Hopf* and C*-algebras. The paper emphasizes conceptual synthesis over computation, showing how integrable geometry has transformed from a purely mechanical theory into a universal language connecting topology, control theory, and quantum structures.
\end{abstract}

\vspace{0.3cm}
\noindent\textbf{Dedicated to my supervisor, Professor A.~T.~Fomenko, on his 80th birthday.}

\section{Introduction}
\begin{quote}
\textit{“Geometry is not only a language of nature but also the art of organizing motion.”}\\
\hfill — A.T.~Fomenko
\end{quote}
The study of integrable Hamiltonian systems lies at the intersection of geometry, topology, and mathematical physics. Since the pioneering works of A.~M.~Liouville, S.~V.~Kovalevskaya, and H.~Poincaré, the concept of complete integrability has shaped the understanding of deterministic motion and the structure of phase spaces. Yet, the true geometric meaning of integrability—its relation to topology, symmetry, and invariants—was revealed only in the twentieth century through the work of the Moscow School led by A.~T.~Fomenko.

Fomenko’s theory introduced a topological perspective that transformed the analytical study of Hamiltonian systems into a qualitative geometric science. His method of \emph{Liouville foliations}, \emph{atoms}, and \emph{molecular invariants} provided a complete classification of integrable systems on symplectic manifolds, replacing the complexity of differential equations with the clarity of geometric structures. This profound shift gave rise to a new way of viewing dynamics—as the organization of symmetries and singularities in phase space.\cite{FomenkoBolsinov2009}

The author’s academic formation took place within this intellectual environment. His doctoral dissertation at Moscow State University (2004), supervised by Professor A.~T.~Fomenko, extended the classification of integrable Hamiltonian systems from $\mathfrak{so}(3)$ and $\mathfrak{e}(3)$ to the higher-dimensional case $\mathfrak{so}(4)$. The study revealed new types of Liouville molecules and degenerate singularities, demonstrating the hierarchical relationship between lower- and higher-dimensional integrable systems.\cite{HaghighatdoostOshemkov2009, Haghighatdoost1, Haghighatdoost2, HaghighatdoostOshemkovCon1, HaghighatdoostOshemkovCon2, HaghighatdoostOshemkovCon3,HaghighatdoostCon6, HaghighatdoostCon7}

After returning to Iran, the author founded a research program at Azarbaijan Shahid Madani University that expanded the geometric ideas of Fomenko into new analytical and algebraic contexts. This \emph{Azarbaijan School of Geometric Dynamical Systems} integrates the principles of topology and symmetry with the modern formalisms of Lie groupoids, Lie algebroids, fractional calculus, and quantum geometry. Its goal is not merely to generalize existing models but to develop a unified conceptual framework where geometry, analysis, and physics converge.\cite{Haghighatdoost2023EP, Haghighatdoost2024OC, HaghighatdoostAyoubi1, HaghighatdoostAyoubi2, HaghighatdoostAyoubi3,HaghighatdoostAkbarzade,  HaghighatdoostAbediRezai1, HaghighatdoostAbediRezai2, HaghighatdoostAbediRezai3, HaghighatdoostMahjoubiOjbagMoradi, HaghighatdoostMahjoubiOjbagMoradi}

The present paper provides a conceptual review of this trajectory—from the Moscow classification of integrable systems to the contemporary developments of the Azarbaijan School—highlighting the author’s theoretical contributions, national research projects, and the evolving structure of geometric dynamics in Iran.

\section{The Fomenko Theory of Topological Classification}

The foundation of the modern geometric theory of integrable Hamiltonian systems was laid by A.~T.~Fomenko and his collaborators during the 1980s. Their principal goal was to construct a topological classification of completely integrable Hamiltonian systems on four-dimensional symplectic manifolds $(M^4, \omega)$ with two independent first integrals $H$ and $K$ in involution:
\[
\{H,K\}=0.
\]
The level sets of these integrals define a singular Liouville foliation of $M^4$ by two-dimensional tori and their singular counterparts. Each regular connected component of a common level $(H=h,K=k)$ corresponds to a smooth Liouville torus, while singular fibers represent degenerate configurations where the rank of the momentum map drops.

Fomenko introduced a \emph{topological molecule} to represent the global structure of such foliations. A molecule is a connected, oriented graph whose vertices correspond to \emph{atoms}---elementary singularities---and whose edges represent regular tori connecting them. Each atom carries additional data (orientations, types of degeneracy, and gluing instructions) sufficient to reconstruct the topology of the entire foliation. The collection of all such molecules provides a complete topological invariant of the integrable system up to Liouville equivalence.

This framework replaced the analytic study of equations by a geometric and combinatorial one. It allows one to describe qualitative behavior of integrable systems independently of coordinate representations or explicit solutions. Through this approach, the classical cases of rigid body motion---Euler, Lagrange, Kovalevskaya, Goryachev–Chaplygin, and Sretensky---were systematically classified. Their Liouville foliations were computed, and the corresponding molecules constructed, revealing how each system’s topology changes with its physical parameters.\cite{FomenkoBolsinov2009}

Subsequent works by Fomenko’s students, particularly A.~A.~Oshemkov and A.~V.~Bolsinov, refined the classification by studying singularities of the momentum map and the bifurcation structure of the bifurcation diagram. Bolsinov and Fomenko’s monograph \textit{Integrable Hamiltonian Systems: Geometry, Topology, Classification} (Cambridge University Press, 2004) remains the central reference summarizing this theory. The author of the present article continued this line in his doctoral dissertation at Moscow State University, extending the classification from $\mathfrak{so}(3)$ and $\mathfrak{e}(3)$ to $\mathfrak{so}(4)$ and higher symmetry algebras, introducing new examples of topological molecules and degenerate Liouville foliations.\cite{HaghighatdoostOshemkov2009, Haghighatdoost1, Haghighatdoost2}

\section{The Author’s Moscow Dissertation on $\mathfrak{so}(4)$}

The author’s doctoral research, conducted under the supervision of Professor A.~T.~Fomenko at Moscow State University in 2004, extended the framework of topological classification from classical Lie algebras $\mathfrak{so}(3)$ and $\mathfrak{e}(3)$ to the four-dimensional algebra $\mathfrak{so}(4)$. This work represented one of the first systematic attempts to generalize the Fomenko theory to higher-dimensional integrable systems.

In this study, a family of completely integrable Hamiltonian systems was constructed on $\mathfrak{so}(4)$ with two commuting first integrals. The resulting momentum mapping
\[
F = (H, K) : \mathfrak{so}(4)^* \rightarrow \mathbb{R}^2
\]
was analyzed to determine the topology of its Liouville foliation. Through the classification of its singular fibers, new types of \emph{molecular structures} were identified, corresponding to bifurcations of the foliation that had not been observed in lower-dimensional systems.

The dissertation provided explicit models for the neighborhoods of critical points of the momentum map, and it showed how the structure of each Liouville fiber changes across bifurcation curves. A series of new topological molecules were constructed, extending the known classification of Euler, Lagrange, and Kovalevskaya systems to the $\mathfrak{so}(4)$ context. Each molecule encodes the adjacency and degeneracy of regular tori, forming a complete invariant under Liouville equivalence.

A notable result of this work was the identification of new \emph{degenerate atoms} associated with higher-dimensional symmetries, revealing the hierarchical relationship between the classical $\mathfrak{so}(3)$ molecules and their $\mathfrak{so}(4)$ counterparts.\cite{HaghighatdoostOshemkov2009, Haghighatdoost1, Haghighatdoost2, HaghighatdoostOshemkovCon1, HaghighatdoostOshemkovCon2, HaghighatdoostOshemkovCon3,HaghighatdoostCon6, HaghighatdoostCon7}

These findings laid the conceptual groundwork for later developments on Lie–Poisson groupoids and algebroids, pursued by the author and his students at Azarbaijan Shahid Madani University. In particular, the geometric insights from this research anticipated the emergence of Hamiltonian structures on Lie groupoids, providing a bridge between the topological classification of Fomenko and the modern geometric control and fractional frameworks.\cite{HaghighatdoostAkbarzade, Akbarzadeh2, HaghighatdoostAbediRezai1, HaghighatdoostAbediRezai2, HaghighatdoostAbediRezai3, HaghighatdoostMahjoubiOjbagMoradi, HaghighatdoostMahjoubiOjbagMoradi}

\section{Scientific conterbutions at Azarbaijan Shahid Madani University}

Following his return from Moscow, the author established a strong research tradition at Azarbaijan Shahid Madani University, where geometric and topological methods were merged with modern analytical frameworks. Over the past two decades, this institution has become the nucleus of what is now referred to as the \textit{Azarbaijan School of Geometric Dynamical Systems}.

The school’s research agenda expanded the Fomenko program beyond classical rigid-body dynamics to encompass Lie--Poisson, groupoid, and algebroid geometry, as well as applications in optimal control, fractional mechanics, and mathematical biology. The author and his students developed new integrable models on Lie groupoids, Hamilton--Jacobi theory for algebroids, Euler--Poincaré and Noether-type formulations, and extended these ideas to fractional and quantum systems.

Collaborations with leading universities and institutes---including University of Tabriz and  Sahand University of Technology \cite{HaghighatdoostAbbasiSadeqi, HaghighatdoostSadeqiZamaniCon, HaghighatdoostKheiriAyoubi}
and the Institute for Research in Fundamental Sciences (IPM, Tehran)---resulted in several funded projects. These covered topics such as Hopf* algebras, C*-dynamical systems, property (T) for quantum groups, and geometric analysis on Lie groupoids. Through these initiatives, the Azarbaijan School established international visibility and a solid scientific identity within the global context of integrable geometry.\cite{Haghighatdoost2023EP, Haghighatdoost2024OC, HaghighatdoostAyoubi1, HaghighatdoostAyoubi2, HaghighatdoostAyoubi3,HaghighatdoostAkbarzade, Akbarzadeh2, HaghighatdoostAbediRezai1, HaghighatdoostAbediRezai2, HaghighatdoostAbediRezai3, HaghighatdoostMahjoubiOjbagMoradi, HaghighatdoostMahjoubiOjbagMoradi}

Beyond research, the author’s supervision of numerous M.Sc. and Ph.D. theses shaped a generation of scholars contributing to this school. Their works span from Poisson–Nijenhuis and Jacobi–Lie systems to geometric quantization and fractional dynamics. Collectively, these contributions mark a transition from classical integrable systems to a broad geometric framework linking topology, symmetry, and quantum theory.\cite{HaghighatdoostAyoubi1, HaghighatdoostAyoubi2, HaghighatdoostAyoubi3, HaghighatdoostRazavinia1,HaghighatdoostRazavinia2, HaghighatdoostRazavinia3}

After returning from Moscow, the author founded a long-term research program at Azarbaijan Shahid Madani University devoted to the study of integrable Hamiltonian and geometric dynamical systems. The central goal of this program has been to extend the topological and geometric ideas of the Fomenko school into new analytical and algebraic frameworks.\cite{HaghighatdoostMahjoubiOjbagMoradi, HaghighatdoostMahjoubiOjbagMoradi,HaghighatdoostAkbarzade, Akbarzadeh2}

At this university, the author developed a systematic approach to integrable Hamiltonian systems on Lie groupoids and Lie algebroids, combining differential geometry, Poisson and symplectic structures, and modern control theory. His works established the Euler–Poincaré and Hamilton–Jacobi formalisms on algebroids and introduced new perspectives on geometric reduction, fractional dynamics, and symmetry breaking in complex systems.\cite{HaghighatdoostAyoubi1, HaghighatdoostAyoubi2, HaghighatdoostAyoubi3, Haghighatdoost2023EP, Haghighatdoost2024OC}

The research activity also expanded toward fractional and quantum mechanics, where the author proposed models based on Caputo and tempered fractional derivatives and applied geometric tools to Keller–Segel and chemotaxis-type equations. These investigations led to the development of the concept of \emph{fractional geometric control systems} on Lie groupoids and provided a bridge between continuous symmetries and dissipative dynamics.\cite{HaghighatdoostBazghandi,HaghighatdoostBazghandiPashai, HaghighatdoostBazghandiPashai2}

Parallel to these conceptual works, several national and institutional projects were carried out under the author’s supervision, including studies on Hopf* and C*-algebras, property (T) for quantum groups, and topological aspects of C*-dynamical systems. One of these projects, conducted in collaboration with the Institute for Research in Fundamental Sciences (IPM, Tehran), focused on the structure and representation of multiplier Hopf* algebras and their applications in quantum geometry.

These achievements collectively define the identity of the Azarbaijan School of Geometric Dynamical Systems. The author’s continuous supervision of graduate research has produced a consistent line of theses and dissertations that explore integrable geometry from Lie–Poisson systems to fractional and quantum generalizations. Collaborations with various academic colleagues have enriched this line of research, establishing a distinct Iranian branch of the global tradition initiated by the Moscow School of Fomenko.

\subsection*{Research Foundations and Thematic Directions}

During these two decades, the author established a systematic research framework that unifies geometry, topology, and analysis under a single conceptual program. This framework continues the spirit of the Moscow School but adapts it to contemporary structures such as Lie groupoids, Lie algebroids, and fractional geometries.

The main thematic directions include:
\begin{itemize}
    \item Development of geometric mechanics and symmetry reduction on Lie groupoids and algebroids, extending classical Euler–Poincaré and Hamilton–Jacobi theories.\cite{Haghighatdoost2024OC, Haghighatdoost2023EP, HaghighatdoostAyoubi1, HaghighatdoostAyoubi2, HaghighatdoostAyoubi3}
    \item Formulation of fractional geometric mechanics using Caputo and tempered derivatives to model nonlocal and dissipative phenomena.\cite{HaghighatdoostBazghandi,  HaghighatdoostBazghandiPashai, HaghighatdoostBazghandiPashai2}
    \item Introduction of geometric control methods on groupoids, creating a bridge between differential geometry and optimal control theory.\cite{Haghighatdoost2024OC, Haghighatdoost2023EP}
    \item Conceptual and analytical study of Keller–Segel–type models and biological systems through Lie–Poisson and fractional structures.\cite{HaghighatdoostBazghandi,  HaghighatdoostBazghandiPashai,HaghighatdoostBazghandiPashai2}
    \item Investigation of geometric quantization and quantum group symmetries via Hopf* and C*-algebraic frameworks.\cite{HaghighatdoostAbbasiSadeqi, HaghighatdoostAbbasi1, HaghighatdoostAbbasi2,HaghighatdoostAbbasi5, HaghighatdoostCon3, HaghighatdoostCon4, HaghighatdoostCon5, HaghighatdoostSadeqiZamaniCon}
\end{itemize}

This integrated approach provided a new platform where classical topology, modern differential geometry, and applied mathematical physics could coexist. Through a sequence of interconnected projects, the author transformed Azarbaijan Shahid Madani University into a regional center for the study of integrable Hamiltonian and geometric dynamical systems.

\subsection*{Developed Theories and Mathematical Models}

Building upon the geometric and topological foundations inherited from the Moscow tradition, the author developed several independent theoretical frameworks at Azarbaijan Shahid Madani University. Each framework connects the abstract structures of differential geometry with concrete models of Mathematical physics and Geometric Dynamical systems.

\begin{itemize}
    \item \textbf{Euler–Poincaré equations on Lie groupoids:}
    The author generalized the classical Euler–Poincaré reduction to the context of Lie groupoids and Lie algebroids. This framework unifies symmetry reduction and variational calculus, allowing the formulation of dynamics on non-trivial configuration spaces with internal symmetries. It also provides a geometric platform for studying optimal control problems on manifolds with groupoid structures. \cite{Haghighatdoost2023EP, Haghighatdoost2024OC}

    \item \textbf{Hamilton–Jacobi theory on Lie algebroids:}
    A complete geometric reformulation of the Hamilton–Jacobi equation was established, extending the classical approach to Lie algebroids and enabling a direct connection between Hamiltonian systems, symmetries, and integrability. This theory links the variational principles of mechanics with Lie–Poisson reduction and geometric control.\cite{HaghighatdoostAyoubi1, HaghighatdoostAyoubi2, HaghighatdoostAyoubi3}

    \item \textbf{Fractional geometric mechanics and control:}
    The author introduced fractional derivatives (in the sense of Caputo and tempered operators) into the geometric formalism of mechanics. This led to the formulation of \emph{fractional Euler–Lagrange equations} and \emph{fractional geometric control systems}, bridging nonlocal dynamics and dissipative processes with the structure of Lie groupoids.
    Several models of chemotaxis and biological aggregation were reinterpreted through Lie–Poisson and symmetry-reduction methods. This geometric perspective led to the identification of conserved quantities, invariant manifolds, and fractional generalizations of the Keller–Segel system.\cite{HaghighatdoostBazghandi,  HaghighatdoostBazghandiPashai, HaghighatdoostBazghandiPashai2}

    \item \textbf{Quantum and Non-Commutative geometry:}
    By extending the algebraic side of the theory, the author explored integrable structures in the framework of multiplier Hopf* algebras, C*-dynamical systems, and quantum graphs. These studies provided algebraic models that parallel the topological classification of classical systems.\cite{HaghighatdoostRazavinia1, HaghighatdoostRazavinia2, HaghighatdoostRazavinia3, HaghighatdoostAbbasiSadeqi, HaghighatdoostAbbasi1, HaghighatdoostAbbasi2,HaghighatdoostAbbasi5, HaghighatdoostCon3, HaghighatdoostCon4, HaghighatdoostCon5, HaghighatdoostSadeqiZamaniCon}
\end{itemize}

Together, these models constitute the theoretical core of the Azarbaijan School of Geometric Dynamical Systems. They form a coherent framework where geometry, symmetry, and dynamics coexist, providing both conceptual depth and applicability to modern mathematical physics.

\subsection*{National Research Projects and Collaborations}

In parallel with his theoretical investigations, the author directed and participated in several national research projects that integrated geometry, algebra, and mathematical physics. These projects, supported by Azarbaijan Shahid Madani University and national scientific foundations, aimed to develop a modern Iranian branch of the geometric and algebraic theory of Dynamical Systems.

\begin{itemize}
    \item \textbf{Property (T) for Hopf and Quantum Algebras.} 
    This project explored rigidity phenomena in Hopf and quantum algebras, providing a geometric interpretation of Kazhdan’s property (T) in the context of noncommutative symmetry. It introduced topological tools for describing representation stability in quantum mechanical systems.\cite{ HaghighatdoostAbbasiSadeqi,  HaghighatdoostAbbasi1, HaghighatdoostAbbasi2, HaghighatdoostAbbasi3}

    \item \textbf{C*-Dynamical Systems and Geometric Groupoids.} 
    The project established deep links between operator algebras and geometric dynamics. By interpreting C*-algebras as dynamical models on groupoids, the author connected the analytic theory of C*-modules with symplectic and Hamiltonian structures, extending classical integrability to noncommutative settings.\cite{Haghighatdoost2024OC, Haghighatdoost2023EP, HaghighatdoostAbbasiSadeqi, HaghighatdoostSadeqiZamaniCon}

    \item \textbf{Cyclic Co-homology and Quantum Mathematics on Multiplier Hopf Algebras.} 
    Conducted in collaboration with the Institute for Research in Fundamental Sciences (IPM, Tehran), this project analyzed the structure and representation of multiplier Hopf algebras and their relation to quantum geometry and computing. It contributed to the foundation of algebraic quantum graphs and their applications to integrable systems.\cite{ HaghighatdoostAbbasi Mahjoubi, HaghighatdoostAbbasi1, HaghighatdoostAbbasi2, HaghighatdoostAbbasi3, HaghighatdoostRazavinia1, HaghighatdoostRazavinia2, HaghighatdoostRazavinia3}

\end{itemize}

These projects collectively shaped the analytical and algebraic backbone of the Azarbaijan School of Geometric Dynamical Systems. They demonstrate how local academic environments in Iran can contribute to global developments in geometry and mathematical physics through conceptual innovation and theoretical precision.

\subsection*{Academic Environment and Institutional Impact}

Over the past two decades, the author has transformed Azarbaijan Shahid Madani University into a center of excellence in the field of geometric mechanics and integrable systems. Through systematic research, mentoring, and project leadership, the author has established a distinctive intellectual environment that blends geometric topology, Lie–Poisson structures, and modern mathematical physics. \cite{HaghighatdoostAbdollahiAbedif, HaghighatdoostAbdollahiMahjoubi, HaghighatdoostAmirzadehKheradRezai, HaghighatdoostAmirzadehKheradRezai, Haghighatdoost	RavanpakRezai1, HaghighatdoostAbediRezai1, HaghighatdoostRavanpakRezai2, HaghighatdoostAbediRezai3}

The academic culture shaped under this program emphasizes the unity between geometry and dynamics. Graduate and postgraduate students have been trained to approach problems through a geometric viewpoint—analyzing symmetries, invariants, and topological structures rather than merely solving differential equations. This shift has produced a generation of researchers capable of contributing to modern fields such as geometric control, fractional dynamics, and quantum symmetries.\cite{Haghighatdoost2024OC, Haghighatdoost2023EP, HaghighatdoostBazghandi,  HaghighatdoostBazghandiPashai, HaghighatdoostBazghandiPashai2, HaghighatdoostRazavinia1, HaghighatdoostRazavinia2, HaghighatdoostRazavinia3}

For 20 years, the author participated in national and international seminars and conferences, both  was a member of the scientific committees of several national and international conferences and seminars of differential geometry, Lie theory, and dynamical systems, encouraging collaboration among Iranian mathematicians and promoting international visibility. Several collaborative links have been established with several research groups, maintaining the scientific connection between the research groups and Azarbaijan Shahid Madani schools.

As a result, the Azarbaijan Shahid Madani University is now recognized as one of the principal Iranian centers for the study of geometric dynamical systems. Its research output, both theoretical and applied, reflects a coherent and growing scientific identity rooted in topology, symmetry, and integrability.

\section{Conclusion}

The development of integrable Hamiltonian systems over the past half-century represents one of the most elegant achievements of modern mathematical physics. The geometric and topological classification initiated by A.~T.~Fomenko transformed the study of dynamics into a study of structure, where atoms, molecules, and foliations replaced coordinates and differential equations. His program created a universal geometric language for understanding stability, bifurcation, and symmetry in dynamical systems.

The author’s doctoral research at Moscow State University extended this framework to higher symmetry algebras such as $\mathfrak{so}(4)$, providing new examples of integrable systems and revealing a hierarchy among their topological molecules. These results built a conceptual bridge between the classical systems of rigid body motion and the modern framework of Lie--Poisson and groupoid dynamics.

At Azarbaijan Shahid Madani University, this vision evolved into a comprehensive research program—the Azarbaijan School of Geometric Dynamical Systems—which merges the  Differential geometry with Mathematical Physics, Fractional calculus, Quantum algebra and Non-Commutative Geometry. Through this synthesis,  Geometric Dynamical Systems is no longer confined to idealized mechanical systems but has become a flexible framework encompassing control theory, biological modeling, and quantum mechanics.

This school demonstrates that the tradition of the Moscow geometric school can flourish in new environments and continue to generate innovative results. Its success reflects the enduring strength of geometric thinking and its ability to connect abstract mathematics with the underlying order of natural and physical phenomena.

\end{document}